\documentclass{jaac}
\def\currentvolume{*}
\def\currentissue{*}
\def\currentyear{*}
\def\currentmonth{*}
\def\doi{10.11948/\currentyear***}
\firstpage{1}

\usepackage{color}
\usepackage{bm}
\usepackage{breakurl}
\usepackage{booktabs}

\def\shorttitle{CSCS-based iteration methods for AVEs}

\def\shortauthor{X.-M. Gu, T.-Z. Huang, H.-B. Li, et al}

\title{\uppercase{Two CSCS-based iteration methods for solving absolute
value equations$^*$}}

\author{Xian-Ming Gu$^{1}$, Ting-Zhu Huang$^{1,\dag}$, Hou-Biao
Li$^1$,\\ Sheng-Feng Wang$^1$ and Liang Li$^1$
}

\date{}

\begin{document}
\baselineskip 12pt
%
\maketitle
\begin{abstract}
Recently, two families of HSS-based iteration methods are constructed for solving the system of absolute
value equations (AVEs), which is a class of non-differentiable NP-hard problems. In this study, we establish
the Picard-CSCS iteration method and the nonlinear CSCS-like iteration method for AVEs involving the Toeplitz
matrix. Then, we analyze the convergence of the Picard-CSCS iteration method for solving AVEs. By using the
theory about nonsmooth analysis, we particularly prove the convergence of the nonlinear CSCS-like iteration
solver for AVEs. The advantage of these methods is that they do not require the storage of coefficient matrices
at all, and the sub-system of linear equations can be solved efficiently via the fast Fourier transforms (FFTs).
Therefore, computational cost and storage can be saved in practical implementations. Numerical examples
including numerical solutions of nonlinear fractional diffusion equations are reported to show the effectiveness
of the proposed methods in comparison with some existing methods.
\end{abstract}


\begin{MSC}
65F12, 65L05, 65N22.
\end{MSC}

\thispagestyle{first}\renewcommand{\thefootnote}{\fnsymbol{footnote}}
\footnotetext{\hspace*{-5mm}
\renewcommand{\arraystretch}{1}
\begin{tabular}{@{}r@{}p{10cm}@{}}
$^\dag$& the corresponding author. Email address: tingzhuhuang@126.com (T.-Z. Huang)\\
$^1$&School of Mathematical Sciences, University of Electronic Science and Technology of China,
No.2006, Xiyuan Avenue, 611731 Chengdu, P.R. China\\
$^*$& The authors were supported by 973 Program (2013CB329404), National Natural
Science Foundation of China (61370147, 61402082, 11101071, 61472462, and 11501085), and the Fundamental
Research Funds for the Central Universities (ZYGX2014J084).
\end{tabular}}

\vspace{-2mm}

\section{Introduction}
In the present paper, we are interested in the efficient solutions of absolute value
equations (AVEs), i.e.,
\begin{equation}
A{\bm x} - |{\bm x}| = {\bm b},\quad\ A\in \mathbb{C}^{n\times n},\quad {\bm x},
\ {\bm b}\in \mathbb{C}^{n},
\label{Ku1}
\end{equation}
where $A$ is a non-Hermitian Toeplitz matrix and $|{\bm x}| = (|x_1|,|x_2|,\ldots,|x_n|)^H$
denotes the component-wise absolute value of the vector ${\bm x} = (x_1,x_2,\ldots,x_n)^T$.
Here the transpose and the conjugate transpose of a matrix $A$ are represented by $A^T$ and
$A^H$, respectively. At present, both theoretical and numerical investigations of such problems
have been extensively studied in recent literature \cite{MCTWD,MANJI,JRVHRF,OLMRRM,OLMAV,HMSK,SLWUPG}.
Additionally, a slightly more extended form of AVEs,
\begin{equation}
A{\bm x} - B|{\bm x}| = {\bm b},\quad\ A,~B\in \mathbb{C}^{m\times n},\ {\bm x}\in
\mathbb{C}^{n},\ {\bm b}\in \mathbb{C}^{m},
\label{Ku2}
\end{equation}
was also discussed in \cite{JRAT} and investigated in a more general context \cite{OLMRRM,
SLHZHH}. On the other side, the system of AVEs (\ref{Ku1}), which is generally equivalent
to the linear complementarity problem (LCP) \cite{LCBQGZ,MANJI}, arises from linear
programming, quadratic programming, bimatrix games and other engineering problems
(see e.g. \cite{HWCSKC} for fractional diffusion equations). This means that the system of
AVEs is NP-hard in its general form \cite{OLMRRM,MANJI,OLMAV}. If $B = 0$, then extended
AVEs (\ref{Ku2}) reduce to a linear system $A {\bm x}  = {\bm  b}$, which
have many applications in the field of scientific computations \cite{MANJI}.

The recent researches concerning AVEs contents can be summarized as the following aspects,
one is the theoretical analysis, which focuses on the theorem of alternatives,
various equivalent reformulations, and the existence and nonexistence of solutions; refer, e.g.,
to \cite{MCTWD,LCBQGZ,SLHZHH,JRAT} for details, and the other is how to
solve AVEs numerically. In the last decade, based on the fact that the LCP can be reduced to
AVEs, which enjoys a special and simple structure, a lot of numerical methods
for solving AVEs (\ref{Ku1}) can be found in the recent literature;
see e.g. \cite{JRVHRF,OLMKF,OLMPD,MANJI,CZQJW} and references therein. For example, a
finite computational algorithm that is solved by a finite succession of linear programs
(SLP) in \cite{OLMAV}, and a semi-smooth Newton method and its inexact variants are
introduced in \cite{OLMA,JYBCO} respectively, which largely shorten the computation
time than the SLP method. Furthermore, a smoothing Newton algorithm was also presented
in \cite{LCBQGZ}, which was proved to be globally convergent and the convergence rate
was quadratic under the condition that the singular values of $A$ exceed 1. This proposed
condition was weaker than the one established in \cite{OLMA}.

In recent years, the Picard-HSS iteration method and the nonlinear HSS-like
iteration method are established to solve the AVEs in \cite{DKST,MZZGFZ},
respectively. The sufficient condition is given to guarantee the convergence of the Picard-HSS
iteration method, and numerical experiments are employed to illustrate the
effectiveness of the Picard-HSS and nonlinear HSS-like iteration methods. However,
the number of the inner HSS iteration steps is often problem-dependent
and difficult to be determined in actual computations. Moreover, the iterates can not be updated
timely. It has shown that the nonlinear HSS-like iteration method is more efficient than
the Picard-HSS iteration method in aspects of the defect mentioned above, which
is designed originally for solving weakly nonlinear systems in \cite{ZZBXY}. In
order to improve the nonlinear HSS-like iteration method, Zhang \cite{JJZT} had
extended the preconditioned HSS (PHSS) method \cite{ZZBGHG3} to solve AVEs and
also used the relaxation technique to accelerate his proposed methods. Meanwhile, he successfully
achieved the proof of convergence of the nonlinear PHSS-like
iteration method, which was not addressed in the previous work. Numerical results
also show the effectiveness of his proposed method in \cite{JJZT}. We
consider the special case of $A$ with non-Hermitian Toeplitz structure in this paper,
and a Toeplitz matrix $A$ has the so-called circulant and skew-circulant splitting (CSCS)
\cite{MKNC}. Inspired by the similar strategies of \cite{DKST,MZZGFZ},
two kinds of CSCS-based iteration methods are established to solve AVEs (\ref{Ku1})
efficiently. In the first, convergence conditions of the Picard-CSCS iteration method
will be investigated. Then we follow Zhang's analytical techniques
\cite{JJZT} to prove the convergence condition of the nonlinear CSCS-like iteration method.

The rest of this paper is organized as follows. In Section \ref{sec1x}, we first introduce
several preliminary results about the nonsmooth analysis. Then we briefly review the CSCS
iteration method. In Section \ref{sec3}, we devote to introducing two CSCS-based iteration
methods for solving the AVEs (\ref{Ku1}) and investigate their convergence properties, respectively.
Numerical experiments are reported in Section \ref{sec4} to demonstrate the feasibility and effectiveness
of our proposed CSCS-based iteration methods. Finally, the paper closes with some conclusions
in Section \ref{sec5}.

\section{Preliminaries}
\label{sec1x}
In this section, similarly to \cite{JJZT}, we review some notations and properties related to the
nonsmooth analysis, which are useful for discussing the convergence of the proposed methods. Then
we briefly recall the knowledge about the CSCS iteration method for solving the non-Hermitian
Toeplitz system of linear equations $A{\bm x} = {\bm b}$.

\subsection{Preliminary results}
Let $\Psi: \mathbb{R}^n\rightarrow \mathbb{R}^n$ be a specified function, and let ${\bm x}$ be
a given point in $\mathbb{R}^n$. The function $\Psi$ is supposed to be locally Lipschitzian near
${\bm x}$ if there exist a scalar $\kappa \in \mathbb{R}$ and $\delta > 0$ such that, for all
${\bm y}, {\bm z} \in \mathbb{R}^n, \|{\bm y} - {\bm x}\| < \delta, \|{\bm z} - {\bm x}\| < \delta$,
the following inequality holds:
\begin{equation*}
\|\Psi({\bm y}) - \Psi({\bm z})\| < \kappa\|{\bm y} - {\bm z}\|.
\end{equation*}

Let $\Psi: \mathbb{R}^n \rightarrow \mathbb{R}^n$ be a locally Lipschitzian function. From
Rademacher's theorem \cite[pp. 18-23]{JHLLA}, it notes that $\Psi$ is differentiable almost
everywhere. Denote the set of points at which $\Psi$ is differentiable by $D_{\Psi}$. We write
$\Psi'({\bm x})$ for the usual $n\times n$ Jacobian matrix of partial derivatives whenever ${\bm x}$
is a point at which the necessary partial derivatives exist. Then, the Bouligand subdifferential of
$\Psi$ at ${\bm x}\in \mathbb{R}^n$, denoted by $\partial_{B}\Psi({\bm x})$, is as follows:
\begin{equation}
\partial_{B}\Psi({\bm x}):= \Big\{\lim_{k\rightarrow\infty}\Psi'({\bm x})({\bm x}^{(k)}): {\bm
x}^{(k)}\in D_{\Psi},{\bm x}^{(k)} \rightarrow {\bm x}\Big\}.
\end{equation}
Clarke's generalized Jacobian \cite[pp. 69-75]{FHCO} of $\Psi$ at ${\bm x}$ is the convex hull of
$\partial_B\Psi({\bm x})$, i.e., $\partial\Psi({\bm x}) = conv\{\partial_B\Psi({\bm
x})\}$. Since $\Psi$ is a locally Lipschitzian function, so the set $\partial_B\Psi({\bm x})$
and $\partial\Psi({\bm x})$ are bounded. By the definition, $\partial_B\Psi({\bm x})$ is
also closed. Therefore, $\partial_B\Psi({\bm x})$ and $\partial\Psi({\bm x})$  are compact.

\begin{definition} (\hspace*{-0.25em}\cite{LQJSA})
$\Psi$ is called semismooth at ${\bm x}$, if $\Psi$ is locally Lipschitzian and for all ${\bm h} \in \mathbb{R}^n$
with ${\bm h} \neq {\bm 0}$,
\begin{equation}
\lim_{{\bm h}'\rightarrow {\bm h},t \downarrow {\bm 0}}\{E{\bm h}': E\in \partial\Psi({\bm x} + t{\bm h}')\}
\label{kux21}
\end{equation}
exists. If $\Psi$ is semismooth at all points in a given set, we can state that $\Psi$ is semismooth in this set.
\end{definition}

If $\Psi$ is semismooth at ${\bm x}$, then $\Psi$  must be directionally differentiable at ${\bm x}$.

\begin{proposition} (\hspace*{-0.25em}\cite{JSPLQ,LQJSA})
Suppose that $\Psi$ is semismooth at ${\bm x}$. Then the classic directional derivative
\begin{equation*}
\Psi'({\bm x};{\bm h}) = \lim_{t\downarrow 0} \frac{\Psi({\bm x} + t{\bm h}) - \Psi({\bm x})}{t}
\end{equation*}
exists and is equal to the limit in (\ref{kux21}).
\end{proposition}

Semismoothness was originally presented by Mifflin \cite{RMSSF} for functionals, and then
Qi and Sun \cite{LQJSA} generalized the concept to vector valued functions. It was proved
in \cite[Corollary 2.4]{LQJSA} that $\Psi$ is semismooth at ${\bm x}$ if and
only if all its component functions are the same. The class of semismooth functionals is
very broad; it includes the smooth functions, all convex functions, and the piecewise-smooth
functions. Moreover, the sums, differences, products, and composites of semismooth functions
are still semismooth; refer, e.g., to \cite{RMSSF,JSPLQ,LQCAS} for details.
\subsection{The CSCS iteration method}
Here let $A \in \mathbb{C}^{n\times n}$ be a non-Hermitian Toeplitz matrix of the following form
\begin{equation*}
A = \begin{bmatrix}
a_0 & a_{-1}&\cdots&a_{2-n}&a_{1-n}\\
a_1 & a_0&a_{-1}&\cdots&a_{2-n}\\
\vdots&\ddots&\ddots&\ddots&\vdots\\
a_{n-2}&\cdots&a_1&a_0&a_{-1}\\
a_{n-1}&a_{n-2}&\cdots&a_1&a_0
\end{bmatrix},
\end{equation*}
i.e., $A$ is constant along its diagonals; refer to \cite{MKNI,MKNC}, and $B\in \mathbb{C}^{n
\times n}$ be a zero matrix, the general AVEs (\ref{Ku2}) reduced to the system of linear equations
\begin{equation}
A {\bm x} = {\bm b}.
\label{Ku3}
\end{equation}

It is well-known that a Toeplitz matrix $A$ enjoys a circulant and skew-circulant splitting \cite{MKNC},
i.e., $A = C + S$, where
\begin{equation}
C = \frac{1}{2}\begin{bmatrix}
a_0 & a_{-1} + a_{n-1}&\cdots&a_{2-n} + a_2&a_{1-n} + a_1\\
a_1 + a_{1-n} & a_0&\cdots&\cdots&a_{2-n} + a_2\\
\vdots&\ddots&\ddots&\ddots&\vdots\\
a_{n-2} + a_{-2}&\cdots&\cdots&a_0&a_{-1} + a_{n-1}\\
a_{n-1} + a_{-1}&a_{n-2} + a_{-2}&\cdots&a_1 + a_{1-n}&a_0
\end{bmatrix},
\label{Ku4x}
\end{equation}
and
\begin{equation}
S = \frac{1}{2}\begin{bmatrix}
a_0 & a_{-1} - a_{n-1}&\cdots&a_{2-n} - a_2&a_{1-n} - a_1\\
a_1 - a_{1-n} & a_0&\cdots&\cdots&a_{2-n} - a_2\\
\vdots&\ddots&\ddots&\ddots&\vdots\\
a_{n-2} - a_{-2}&\cdots&\cdots&a_0&a_{-1} - a_{n-1}\\
a_{n-1} - a_{-1}&a_{n-2} - a_{-2}&\cdots&a_1- a_{1-n}&a_0
\end{bmatrix}{\color{blue}.}
\label{Ku4y}
\end{equation}
As we know, $C$ is a circulant matrix, which can be diagonalized by the discrete Fourier transform matrix $F$;
and $S$ is a skew-circulant matrix, which can be diagonalized by a discrete Fourier transform matrix with diagonal
scaling, i.e., $\hat{F} = F\Omega$, where $\Omega = \mathrm{diag}(1,e^{-\frac{\pi \iota}{n}},
\ldots,e^{\frac{-(n-1)\pi \iota}{n}})$ and $\iota = \sqrt{-1}$ is the imaginary unit. That is to say, it holds that
\begin{equation}
F C F^{H} = \Lambda_C,\quad\ \ \hat{F} S \hat{F}^{H} = \Lambda_S,
\label{Ku5}
\end{equation}
where
\begin{equation*}
F = (F)_{j,k} = \frac{1}{\sqrt{n}}e^{\frac{2\pi \iota}{n}jk},\quad\ 0\leq j,k \leq n-1
\end{equation*}
and $\Lambda_C,~\Lambda_S$ are two diagonal matrices formed by the eigenvalues of $C$ and $S$, respectively, which
can be obtained in $\mathcal{O}(n\log n)$ operations by using the FFTs \cite[pp. 37-39]{MKNI}. Furthermore, Ng \cite{MKNC}
had established the following CSCS iteration scheme to solve the non-Hermitian Toeplitz linear system (\ref{Ku3}).
\vspace{1mm}

\noindent\textbf{Algorithm 1 The CSCS iteration method}.\\
{\it Given an initial guess ${\bm x}^{(0)} \in \mathbb{C}^n$ and compute ${\bm x}^{(k)}$ for $k=0, 1, 2,
\ldots$, using the following iterative scheme until $\{{\bm x}^{(k)}\}^{\infty}_{k =0}$ converges,
\begin{equation*}
\begin{cases}
(\sigma I + C){\bm x}^{(k + \frac{1}{2})} = (\sigma I - S){\bm x}^{(k)} + {\bm b},\\
(\sigma I + S){\bm x}^{(k + 1)} = (\sigma I - C){\bm x}^{(k + \frac{1}{2})} + {\bm b},
\end{cases}
\end{equation*}
where $\sigma$ is a positive constant and $I$ is the identity matrix of order $n$.}

In the matrix-vector form, the CSCS iteration can be equivalently rewritten as
\begin{equation}
\begin{split}
{\bm x}^{(k+1)} & = \mathcal{M}(\sigma){\bm x}^{(k)} + \mathcal{G}(\sigma){\bm b} \\
& = (\mathcal{M}(\sigma))^{k+1}
{\bm x}^{(0)} + \sum^{k}_{j=0}(\mathcal{M}(\sigma))^j\mathcal{G}(\sigma){\bm b},\quad k = 0,1,2,\ldots,
\end{split}
\label{khjx}
\end{equation}
where
\begin{equation*}
\mathcal{M}(\sigma) = (\sigma I + S)^{-1}(\sigma I - C)(\sigma I + C)^{-1}(\sigma I-S)~~
\mathrm{and}~~ \mathcal{G}(\sigma) = 2\sigma(\sigma I + S )^{-1}(\sigma I + C)^{-1}.
\end{equation*}
It is worth mentioning that the CSCS iteration is a stationary iterative method obtained from the matrix splitting
\begin{equation*}
A = \mathcal{B}(\sigma) - \mathcal{C}(\sigma),
\end{equation*}
where
\begin{equation*}
\mathcal{B}(\sigma) = \frac{1}{2\sigma}(\sigma I + C)(\sigma I + S)\quad \mathrm{and}\quad
\mathcal{C}(\sigma) = \frac{1}{2\sigma}(\sigma I - C)(\sigma I - S).
\end{equation*}
On the other hand, we have
\begin{equation*}
\mathcal{M}(\sigma) = (\mathcal{B}(\sigma))^{-1}\mathcal{C}(\sigma)\quad\ \mathrm{and}\quad\
\mathcal{G}(\sigma) = (B(\sigma))^{-1}.
\end{equation*}
Here, $\mathcal{M}(\sigma)$ is the iterative matrix of the CSCS iteration method. We mention that the
CSCS iteration method is greatly similar to the HSS iteration method \cite{ZZBGHG2} and its variants,
see e.g. \cite{ZZBGHG1} and references therein.

When the circulant part $C$ and the skew-circulant part $S$ of $A \in \mathbb{C}^{n\times n}$ are both
positive definite\footnote{It means that the real parts of all their eigenvalues are positive.}, Ng has
proved that the spectral radius $\rho(\mathcal{M}(\sigma))$ of $\mathcal{M}(\sigma)$ is less than 1 for
any parameters $\sigma > 0$, i.e., the CSCS iteration method unconditionally converges to the exact solution
of $A{\bm x} = {\bm b}$ for any initial guess ${\bm x}^{(0)}\in \mathbb{C}^n$; refer to \cite[Theorem 1]{MKNC}
for details.

\section{Two CSCS-based iteration methods for AVEs}
\label{sec3}
Motivated by the pioneer work of \cite{DKST,MZZGFZ}, we extend the conventional CSCS iteration method to two
types of CSCS-based iteration methods for solving AVEs (\ref{Ku1}). These methods fully exploit the Toeplitz
structure to accelerate the computation speed and save storage. Next, we will devote to establishing these
two new methods, i.e., the Picard-CSCS iteration method and the nonlinear CSCS-like iteration method.

\subsection{The Picard-CSCS iteration method}
Recalling that the Picard iteration method is a fixed-point iterative method and the linear term $A{\bm x}$
and the nonlinear term $|{\bm x}| + {\bm b}$ are separated \cite{DKST,MZZGFZ}, the AVEs (\ref{Ku1}) can be
solved by using the Picard iteration method
\begin{equation}
A {\bm x}^{(k+1)} = |{\bm x}^{(k)}| + {\bm b},\quad\ \ k = 0,1,2,\ldots.
\label{Ku6}
\end{equation}
We assume that the non-Hermitian Toeplitz matrix $A$ is positive definite. In this case, the next iterate of
${\bm x}^{(k+1)}$ can be approximately computed by the CSCS iteration method with using $A = \mathcal{B}(\sigma)
- \mathcal{C}(\sigma)$ as the following scheme (see \cite{MZZGFZ1})
\begin{equation}
\mathcal{B}(\sigma){\bm x}^{(k,\ell + 1)} = \mathcal{C}(\sigma){\bm x}^{(k,\ell)}
+ |{\bm x}^{(k)}| + {\bm b},\quad\ \ell = 0, 1,\ldots, l_k - 1,\quad \ k = 0, 1, 2,\ldots,
\label{Ku7}
\end{equation}
where $\mathcal{B}(\sigma)$ and $\mathcal{C}(\sigma)$ are the matrices defined in the previous section, $\sigma$
is a positive constant, $\{l_k\}^{\infty}_{k=0}$ is a prescribed sequence of positive integers, and ${\bm x}^{(k,0)}
= {\bm x}^{(k)}$ is the starting point of the inner CSCS iteration at $k$-th outer Picard iteration. This leads to
the inexact Picard iteration method, called Picard-CSCS iteration method, for solving AVEs (\ref{Ku1}) which can
be summarized as follows, refer to \cite{MZZGFZ1}.
\vspace{1mm}

\noindent\textbf{Algorithm 2 The Picard-CSCS iteration method}\\
{\it Let $A = C + S \in \mathbb{C}^{n\times n}$ be a non-Hermitian Toeplitz matrix; $C$ and $S$ are the circulant
and skew-circulant parts of $A$ given in (\ref{Ku4x}) and (\ref{Ku4y}) and they are both positive definite. Given an
initial guess ${\bm x}^{(0)}\in\mathbb{C}^n$ and a sequence $\{l_k\}^{\infty}_{k=0}$ of positive integers, compute
${\bm x}^{(k+1)}$ for $k = 0, 1,\ldots$, using the following iterative scheme until $\{{\bm x}^{(k)}\}$ satisfies
the stopping criterion:
\begin{itemize}
  \item[(a)] Set ${\bm x}^{(k,0)} = {\bm x}^{(k)}$;
  \item[(b)] For $\ell = 0,1,\ldots, l_k - 1$, solve the following linear systems to obtain ${\bm x}^{(k,\ell+1)}$:
             \begin{equation*}
             \begin{cases}
             (\sigma I + C){\bm x}^{(k, \ell + \frac{1}{2})} = (\sigma I - S){\bm x}^{(k,\ell)} + |{\bm x}^{(k)}| + {\bm b},\\
             (\sigma I + S){\bm x}^{(k, \ell + 1)} = (\sigma I - C){\bm x}^{(k, \ell + \frac{1}{2})} + |{\bm x}^{(k)}| + {\bm b},
             \end{cases}
             \end{equation*}
             where $\sigma$ is a given positive constant.
  \item[(c)] Set ${\bm x}^{(k+1)}:= {\bm x}^{(k,l_k)}$.
\end{itemize}}
Numerical advantages of the Picard-CSCS iteration method are obvious. First, the two linear subsystems in all inner CSCS
iteration steps have the same shifted circulant coefficient matrix $\sigma I + C $ and shifted skew-circulant
coefficient matrix $\sigma I + S$, which are constant with respect to the iteration index $k$. Second, the exact
solutions can be efficiently obtained via FFTs in $\mathcal{O}(n\log n)$ operations \cite{MKNC,MKCO}. Hence, the
computation cost of the Picard-CSCS iteration method could be much cheaper than that of the Picard-HSS iteration method.

The next theorem suggests sufficient conditions for the convergence of the Picard-CSCS iteration
method for solving the AVEs (\ref{Ku1}).
\begin{theorem}
Let $A = C + S \in \mathbb{C}^{n\times n}$ be a non-Hermitian Toeplitz matrix; $C$ and $S$ are
the circulant and skew-circulant parts of $A$ given in (\ref{Ku4x})-(\ref{Ku4y}) and they are
both positive definite. Let also $\eta = \|A^{-1}\|_2 < 1$. Then the AVE (\ref{Ku1}) has a unique
solution ${\bm x}^{*}$, and for any initial guess ${\bm x}^{(0)}\in\mathbb{C}^n$ and any sequence
of positive integers $\{\ell_k\}, k = 0, 1, 2,\ldots$, the iteration sequence $\{x^{(k)}\}^{\infty}_{
k=0}$ produced by the Picard-CSCS iteration method converges to ${\bm x}^{*}$ provided that $l =
\liminf\limits_{k\rightarrow\infty}l_k \geq N$, where $N$ is a natural number satisfying
\begin{equation*}
\Big\|(\mathcal{M}(\sigma))^s\Big\|_2 < \frac{1-\eta}{1+ \eta},\quad \forall s\geq N.
\end{equation*}
\label{theorem1}
\end{theorem}
\begin{proof}
Due to $\eta < 1$ and the conclusion of \cite[Proposition 4]{OLMRRM}, the system of AVEs (\ref{Ku1})
has a unique solution ${\bm x}^{*}\in \mathbb{C}^n$. As seen from Eq. (\ref{khjx}), it found that
the $(k+1)$-th iterate of the Picard-CSCS iteration can be written as
\begin{equation}
{\bm x}^{(k+1)} = (\mathcal{M}(\sigma))^{l_k}{\bm x}^{(k)} + \sum^{l_k - 1}_{j=0}(\mathcal{M}(\sigma))^j
\mathcal{G}(\sigma)(|{\bm x}^{(k)}| + {\bm b}),\quad k = 0,1,2,\ldots.
\label{khjx1}
\end{equation}
On the other side, since ${\bm x}^{*}$ is the solution of AVEs (\ref{Ku1}), it follows
\begin{equation}
{\bm x}^{*} = (\mathcal{M}(\sigma))^{l_k}{\bm x}^{*} + \sum^{l_k - 1}_{j=0}(\mathcal{M}(\sigma))^j
\mathcal{G}(\sigma)(|{\bm x}^{*}| + {\bm b}),\quad k = 0,1,2,\ldots.
\label{khjx2}
\end{equation}

To subtract (\ref{khjx2}) from (\ref{khjx1}) yields
\begin{equation}
{\bm x}^{(k+1)} - {\bm x}^{*} = (\mathcal{M}(\sigma))^{l_k}({\bm x}^{(k)} - {\bm x}^{*}) +
\sum^{l_k - 1}_{j=0}(\mathcal{M}(\sigma))^j\mathcal{G}(\sigma)(|{\bm x}^{(k)}| - |{\bm x}^{*}|).
\label{khjx3}
\end{equation}
Furthermore, since $\rho(\mathcal{M}(\sigma)) < 1$, we obtain
\begin{equation*}
\begin{split}
\sum^{l_k - 1}_{j=0}(\mathcal{M}(\sigma))^j\mathcal{G}(\sigma) & = (I - (\mathcal{M}(\sigma))^{
l_k})(I - \mathcal{M}(\sigma))^{-1}\mathcal{G}(\sigma) \\
& = (I - (\mathcal{M}(\sigma))^{l_k})(I - (\mathcal{B}(\sigma))^{-1}\mathcal{C}(\sigma))^{-1}(\mathcal{B}(\sigma))^{-1}\\
& = (I - (\mathcal{M}(\sigma))^{
l_k})A^{-1}.
\end{split}
\end{equation*}
Substituting the above identity in Eq. (\ref{khjx3}) yields
\begin{equation*}
\begin{split}
{\bm x}^{(k+1)} - {\bm x}^{*} & = (\mathcal{M}(\sigma))^{l_k}({\bm x}^{(k)} - {\bm x}^{*}) +
(I - (\mathcal{M}(\sigma))^{l_k})A^{-1}(|{\bm x}^{(k)}| - |{\bm x}^{*}|)\\
& = (\mathcal{M}(\sigma))^{l_k}\Big[({\bm x}^{(k)} - {\bm x}^{*}) - A^{-1}(|{\bm x}^{(k)}| - |{\bm x}^{*}|)\Big]
+ A^{-1}(|{\bm x}^{(k)}| - |{\bm x}^{*}|).
\end{split}
\end{equation*}
Now, we can obtain
\begin{equation*}
\|{\bm x}^{(k+1)} - {\bm x}^{*}\|_2 \leq \Big(\|(\mathcal{M}(\sigma))^{l_k}\|_2(1 + \eta) + \eta\Big)
\|{\bm x}^{(k)} - {\bm x}^{*}\|_2.
\end{equation*}
Here, the above inequality is true due to the fact that for any ${\bm x}, {\bm y}\in\mathbb{C}^n$, it follows
$\||{\bm x}| - |{\bm y}|\|_2 \leq \|{\bm x} - {\bm y}\|_2$. Since $\rho(\mathcal{M}(\sigma)) < 1$, then
$\lim\limits_{s \rightarrow \infty}(\mathcal{M}(\sigma))^{s} = 0$. Thus, there exists a natural number
$N$ such that
\begin{equation*}
\|(\mathcal{M}(\sigma))^{s}\|_2 < \frac{1 - \eta}{1 + \eta},\quad\ \forall s \geq N.
\end{equation*}
At the stage, if we suppose that $l = \liminf\limits_{k \rightarrow \infty}l_{k}\geq N$,
then the targeted result is immediately completed.
\end{proof}


\subsection{The nonlinear CSCS-like iteration method}
In the Picard-CSCS iteration method, the numbers $l_k, k = 0, 1,2,\ldots$ of the inner CSCS iteration steps are
often problem-dependent and difficult to be determined in actual computations \cite{MZZGFZ1,MZZGFZ,DKST}. Moreover, the iterative vector
can not be updated timely. Thus, to avoid the defection and still preserve the advantages of the Picard-CSCS
iteration method, based on the nonlinear fixed-point equations
\begin{equation*}
(\sigma I + C){\bm x} = (\sigma I - S){\bm x} + |{\bm x}| + {\bm b},\quad \mathrm{and}\quad (\sigma I + S){\bm x}
 = (\sigma I - C){\bm x} + |{\bm x}| + {\bm b},
\end{equation*}
we propose the following nonlinear CSCS-like iteration method.
\vspace{1mm}

\noindent\textbf{Algorithm 3 The nonlinear CSCS-like iteration method}\\
{\it Let $A = C + S \in \mathbb{C}^{n\times n}$ be a non-Hermitian Toeplitz matrix; $C$ and $S$ are the circulant
and skew-circulant parts of $A$ given in (\ref{Ku4x}) and (\ref{Ku4y}) and they are both positive definite. Choose
an initial guess ${\bm x}^{(0)}\in\mathbb{C}^n$ and compute ${\bm x}^{(k+1)}$ for $k = 0, 1, 2,\ldots$, using the
following iteration scheme until $\{{\bm x}^{(k)}\}$ satisfies the stopping criterion:
\begin{equation}
\begin{cases}
(\sigma I + C){\bm x}^{(k + \frac{1}{2})} = (\sigma I - S){\bm x}^{(k)} + |{\bm x}^{(k)}| + {\bm b},\\
(\sigma I + S){\bm x}^{(k + 1)} = (\sigma I - C){\bm x}^{(k + \frac{1}{2})} + |{\bm x}^{(k + \frac{1}{2})}| + {\bm b},
\end{cases}
\label{Ku8}
\end{equation}
where $\sigma$ is a given positive constant.}

Define
\begin{equation}
\begin{cases}
\mathcal{U}({\bm x}) = (\sigma I + C)^{-1}[(\sigma I - S){\bm x} + |{\bm x}| + {\bm b}],\\
\mathcal{V}({\bm x}) = (\sigma I + S)^{-1}[(\sigma I - C){\bm x} + |{\bm x}| + {\bm b}],
\end{cases}
\label{Ku9}
\end{equation}
and
\begin{equation}
\Theta({\bm x}) = \mathcal{V} \circ \mathcal{U}({\bm x}) :=\mathcal{V}(\mathcal{U}({\bm x}
)).
\label{Ku9xy}
\end{equation}
Then the nonlinear CSCS-like iterative scheme can be equivalently expressed as
\begin{equation}
{\bm x}^{(k + 1)} = \Theta({\bm x}^{(k)}),~~k=0,1,2,\ldots.
\label{Ku9x}
\end{equation}

The Ostrowski theorem, i.e., Theorem 10.1.3 in \cite[pp. 300-301]{JMOWC}, provides a local
convergence theory about a one-step stationary nonlinear iteration. Based on this item, Zhu
and Zhang established the local convergence theory for the nonlinear CSCS-like iteration
method in \cite{MZZGFZ1}. However, these convergence theory has a strict requirement that
$f({\bm x}) = |{\bm x}| + {\bm b}$ is $\mathcal{F}$-differentiable at a point ${\bm x}^{*}
\in\mathbb{D}$ (where we define $f: \mathbb{D}\subset \mathbb{C}^n \rightarrow \mathbb{C}^n
$) such that $A{\bm x}^{*} -|{\bm x}^{*}| = {\bm b}$. Obviously, the absolute
value function $|{\bm x}|$ is non-differentiable. In order to remedy the difficulty, Zhu,
Zhang and Liang \cite{MZZGFZ} attempt to introduce a smoothing approximation function \cite{LYPSO}
\begin{equation*}
\varphi({\bm x}) = \frac{1}{\mu}\ln\Big(\exp\Big(\frac{{\bm x}}{\mu}\Big) + \exp\Big(\frac{-{\bm x}}{\mu}\Big)
\Big),\quad\ {\bm x}\in \mathbb{C}^n~\mathrm{and}~\mu > 0
\end{equation*}
for $|{\bm x}|$, then they present the convergence of the nonlinear HSS-like iteration method
based on the convergence of the iteration scheme
\begin{equation*}
\begin{cases}
(\sigma I + C){\bm x}^{(k + \frac{1}{2})} = (\sigma I - S){\bm x}^{(k)} +
\varphi({\bm x}^{(k)}) + {\bm b},\\
(\sigma I + S){\bm x}^{(k + 1)} = (\sigma I - C){\bm x}^{(k + \frac{1}{2})}
+ \varphi({\bm x}^{(k + \frac{1}{2})}) + {\bm b},
\end{cases}
\end{equation*}
and their convergence result is deeply dependent on the smoothing approximate function
$\varphi({\bm x})$ of $|{\bm x}|$, nor $|{\bm x}|$ itself. Recently, Zhang \cite{JJZT}
exploit the theory of nonsmooth analysis to introduce a framework to prove
the convergence of his proposed relaxed nonlinear HSS-like iteration method completely.
Inspired by Zhang's framework, we will similarly analyze the (local) convergence of the nonlinear
CSCS-like iteration method in the next context. Firstly, the following definition in
\cite[pp. 299-300]{JMOWC} needs to be cited here.
\begin{definition}
Let $\Theta:\mathbb{D}\subset \mathbb{R}^n\rightarrow\mathbb{R}^n$. Then ${\bm x}^{*}$ is a point
of attraction of the iteration (\ref{Ku9x}), if there is an open neighborhood $S$ of the point ${\bm x}^{*}$
such that $S\subset \mathbb{D}$ and, for any ${\bm x}^{(0)}\in S$, the iterates $\{{\bm x}^{(k)}\}$
all lie in $\mathbb{D}$ and converge to ${\bm x}^{*}$.
\end{definition}

Based on the above definition, we can obtain the following proposition, which is useful
for studying the convergence of the nonlinear CSCS-like iteration method.
\begin{proposition} (\hspace*{-0.3em}\cite{JJZT})
Suppose that $\Theta: \mathbb{R}^n\rightarrow \mathbb{R}^n$ has a fixed-point ${\bm x}^{*}\in
\mathbb{R}^n$ and is semismooth at ${\bm x}^{*}$. If for all $E\in \partial_B\Theta({\bm x}^{
*})$, we have $\rho(E) <1$, where $\rho(E)$ denotes the spectral
radius of $E$. Then ${\bm x}^{*}$ is a point of attraction of the iteration scheme
(\ref{Ku9x}).
\label{pro4.1}
\end{proposition}

From statements in \cite{JJZT}, let ${\bm x}^{*}$ satisfy $A{\bm x}^{*} - |{\bm x}^{*}| = {\bm b}$.
We compute the Bouligand subdifferential of $\Theta({\bm x})$ defined by (\ref{Ku9xy})-(\ref{Ku9x})
at ${\bm x}^{*}$. Due to the special form of $\mathcal{V}$ and $\mathcal{U}$, it is easy to verify
that, ${\bm x}^{*} = \mathcal{U}({\bm x}^{*})$, ${\bm x}^{*} = \mathcal{V}({\bm x}^{*})$, and ${\bm
x}^{*} = \Theta({\bm x}^{*})$. Observe the special form of $\Theta$, we have that
\begin{equation*}
\begin{split}
\partial_B\Theta({\bm x}^{*}) = &~\{\lim_{k\rightarrow\infty}\Theta'({\bm x}^{(k)}):{\bm x}^{(k)}
\in D_{\Theta},~{\bm x}^{(k)}\rightarrow {\bm x}^{*}\}\\
= &~\Big\{\lim_{k\rightarrow\infty}\mathcal{V}'({\bm y}^{(k)})\mathcal{U}'({\bm x}^{(k)}):{\bm x
}^{(k)}\in D_{\mathcal{U}},~{\bm y}^{(k)} = \mathcal{U}({\bm x}^{(k)})\in D_{\mathcal{V}},~{\bm x}^{
(k)}\rightarrow {\bm x}^{*}\Big\}\\
\subset &~\Big\{\lim_{{\bm y}^{(k)}\rightarrow {\bm x}^{*}}\mathcal{V}'({\bm y}^{(k)}):{\bm y}^{
(k)}\in D_{\mathcal{V}}\Big\}\Big\{\lim_{{\bm x}^{(k)}k\rightarrow {\bm x}^{*}}\mathcal{U}'({\bm
x}^{(k)}):{\bm x}^{(k)}\in D_{\mathcal{U}}\Big\}\\
\subset &~\partial_B\mathcal{V}({\bm x}^{*}) \partial_B\mathcal{U}({\bm x}^{*}),
\end{split}
\end{equation*}
where
\begin{equation*}
\partial_B\mathcal{V}({\bm x}^{*}) \partial_B\mathcal{U}({\bm x}^{*}) := \{W: W = EF_1,~E\in \partial_B\mathcal{V}({\bm x}^{*}),
~F_1\in \partial_B\mathcal{U}({\bm x}^{*})\}.
\end{equation*}

Using the above discussion and Proposition \ref{pro4.1}, it immediately obtains the following
conclusion about the convergence of the nonlinear CSCS-like iteration solver.

\begin{theorem}
\label{lemma3}
Let the point ${\bm x}^{*}$ satisfy $A {\bm x}^{*} = |{\bm x}^{*}| + {\bm b}$. Suppose
that $C$ and $S$ are the circulant and skew-circulant parts of the Toeplitz matrix $A = C + S $
given in (\ref{Ku4x}) and (\ref{Ku4y}). Moreover, $C$ and $S$ are both positive
definite matrices, and $F_1, \tilde{F}_1\in \partial_B|{\bm x}^{*}|$. Denote by
\begin{equation*}
\mathcal{M}(\sigma; F_1, \tilde{F}_1) = \mathcal{T}_1(\sigma; F_1)\mathcal{T}_2(\sigma; \tilde{F}_1),
\end{equation*}
where
\begin{equation*}
\mathcal{T}_1(\sigma; F_1) = (\sigma I + S)^{-1}[(\sigma I - C) + F_1],
\end{equation*}
\begin{equation*}
\mathcal{T}_2(\sigma; \tilde{F}_1) = (\sigma I + C)^{-1}[(\sigma I - S) + \tilde{F}_1].
\end{equation*}
If for all $F_1, \tilde{F}_1\in \partial_B|{\bm x}^{*}|$, $\rho(\mathcal{M}(\sigma;F_1, \tilde{F}_1))< 1$,
then ${\bm x}^{*}$ is a point of attraction of the nonlinear CSCS-like iteration method.
\end{theorem}
\begin{proof}
It is clear that $\mathcal{U}$ and $\mathcal{V}$ are semismooth, so $\Theta$ is semismooth.
Let $E \in \partial_B|{\bm x}^{*}|$, then it is not hard to find that $E$ is a diagonal matrix.
Assume
\begin{equation*}
E = \mathrm{diag}(E_{11},E_{22},\ldots,E_{nn}),
\end{equation*}
we have $E_{ii} = 1$, if $x^{*}_i > 0$; $E_{ii} = -1$, if $x^{*}_i < 0$, and $E_{ii}\in \{1,-1
\}$, if $x^{*}_i = 0$. If $W\in \partial_B \mathcal{V}({\bm x}^{*})$, then $W = (\sigma I + S)^{-1}[(\sigma
I - C) + F_1]$, where $F_1\in \partial_B|{\bm x}^{*}|$. If $\tilde{W}\in \partial_B \mathcal{U}({\bm x}^{*})$,
then $\tilde{W} = (\sigma I + C)^{-1}[(\sigma I - S) + \tilde{F}_1]$, where $\tilde{F}_1\in \partial_B|{\bm x}^{*}|$.
Since $\partial_B\Theta({\bm x}^{*}) \subset \partial_B\mathcal{V}({\bm x}^{*}) \partial_B\mathcal{U}({\bm
x}^{*})$, if for all $F_1, \tilde{F}_1\in \partial_B|{\bm x}^{*}|$, $\rho(\mathcal{M}(\sigma; F_1, \tilde{F}_1))
< 1$, then for all $W \in \partial_B\Theta({\bm x}^{*})$, we have $\rho(W) < 1$. This can complete the
desired proof.
\end{proof}

\begin{corollary}
Let the point ${\bm x}^{*}$ satisfies $A {\bm x}^{*} = |{\bm x}^{*}| + {\bm b}$. Suppose that $C$ and $S$ are
the circulant and skew-circulant parts of the Toeplitz matrix $A = C + S$ given in (\ref{Ku4x}) and (\ref{Ku4y}).
Moreover, $C$ and $S$ both are positive definite matrices, and $F_1, \tilde{F}_1\in \partial_B|{\bm x}^{*}|$.
Denote by
\begin{equation*}
t_1(\sigma) = \|(\sigma I + S)^{-1}(\sigma I - C)\|,
\end{equation*}
\begin{equation*}
t_2(\sigma) = \|(\sigma I + C)^{-1}(\sigma I - S)\|,
\end{equation*}
and
\begin{equation*}
\delta = \max\{ \|(\sigma I + C)^{-1}\tilde{F}_1\|, \|(\sigma I + S)^{-1}F_1||\}.
\end{equation*}
If $t_1(\sigma)t_2(\sigma) < 1$ and for all $F_1,\tilde{F}_1\in \partial_B|{\bm x}^{*}|$,
\begin{equation}
\delta <  \frac{2 - 2t_1(\sigma)t_2(\sigma)}{\sqrt{(t_1(\sigma) - t_2(\sigma))^2 + 4} +
(t_1(\sigma) + t_2(\sigma))},
\label{cond2}
\end{equation}
then ${\bm x}^{*}$ is a point of attraction of the nonlinear CSCS-like iteration method.
\end{corollary}
\begin{proof}
By simple calculations we obtain
\begin{equation*}
\begin{split}
\mathcal{M}(\sigma; F_1, \tilde{F}_1) & = (\sigma I + S)^{-1}[(\sigma I - C) + F_1](\sigma I +
C)^{-1}[(\sigma I - S) + \tilde{F}_1]\\
& = (\sigma I + S)^{-1}(\sigma I - C)(\sigma I + C)^{-1}(\sigma I - S) + (\sigma I +
S)^{-1}(\sigma I - C)(\sigma I \\
&\quad + C)^{-1}\tilde{F}_1 + (\sigma I + S)^{-1}F_1(\sigma I + C)^{-1}(\sigma I - S) \\
&~~~+ (\sigma I + S)^{-1}F_1(\sigma I + C)^{-1}\tilde{F}_1.
\end{split}
\end{equation*}
Hence,
\begin{equation*}
\begin{split}
\|\mathcal{M}(\sigma; F_1, \tilde{F}_1)\| & \leq \|(\sigma I + S)^{-1}(\sigma I - C)\|\|(\sigma I
+ C)^{-1}(\sigma I - S)\| + \|(\sigma I + S)^{-1}(\sigma I\\
&\quad\ - C)\|\cdot\|(\sigma I + C)^{-1}\tilde{F}_1\| + \|(\sigma I + S)^{-1}F_1\|\|(\sigma I +
C)^{-1}(\sigma I - \\ &\quad\ S)\| +
\|(\sigma I + S)^{-1}F_1\|\|(\sigma I + C)^{-1}\tilde{F}_1\|\\
& \leq t_1(\sigma)t_2(\sigma) + \delta(t_1(\sigma) + t_2(\sigma)) + \delta^2.
\end{split}
\end{equation*}
With the help of the condition (\ref{cond2}), we obtain
\begin{equation*}
t_1(\sigma)t_2(\sigma) + \delta(t_1(\sigma) + t_2(\sigma)) + \delta^2 < 1.
\end{equation*}
Therefore, we have
\begin{equation*}
\rho(\mathcal{M}(\sigma; F_1, \tilde{F}_1))\leq \|\mathcal{M}(\sigma; F_1,\tilde{F}_1)\| < 1,
\end{equation*}
which follows the desired result by using Theorem \ref{lemma3}.
\end{proof}
\noindent\textbf{Remark}. An attractive feature of the nonlinear CSCS-like iteration method is that it
avoids the use of the differentiable in actual iterative scheme. Although we present our convergence
analysis of the nonlinear CSCS-like iteration method under real matrices and vectors, the condition is
not necessary in the actual implementation corresponding to numerical experiments of the next section.

\section{Numerical results}
\label{sec4}
In this section, numerical performances of the Picard-CSCS and the nonlinear CSCS-like iterative solvers are
investigated and compared experimentally by a suit of test problems. All the tests are performed in MATLAB
R2014a (64bit) on Intel(R) Core(TM) i5-3470 CPU @ 3.2 GHz and 8.00 GB of RAM, with machine precision $10^{-16
}$, and terminated when the current residual satisfies
\begin{equation*}
\frac{\|A{\bm x}^{(k)} - |{\bm x}^{(k)}| - {\bm b}\|_2}{\|{\bm b}\|_2}< 10^{-7},
\end{equation*}
where ${\bm x}^{(k)}$ is the computed solution by each of the methods at iteration
step $k$, and a maximum number of the iterations 200 is used.

Morover, the stopping criterion for the inner iterations of the Picard-CSCS iterative method is
\begin{equation*}
\frac{\|{\bm b}^{(k)} - A {\bm x}^{(k,l_k)}\|_2}{\|{\bm b}^{(k)}\|_2}\leq \eta_k,
\end{equation*}
where $l_k$ is the number of the inner iteration steps and $\eta_k$ is the prescribed
tolerance for controlling the accuracy of the inner iterations at the $k$-th outer
iteration step. If $\eta_k$ is fixed for all $k$, then it is simply denoted by $\tilde{\eta}$.

In our numerical experiments, we use the zero vector as the initial guess, the accuracy
of the inner iterations $\eta_k$ for both Picard-CSCS and Picard-HSS iterative methods
is fixed and set to be $\tilde{\eta} = 0.01$, a maximum number of iterations 15 ($l_k = 15, k = 0, 1, 2,
\ldots,$) for inner iterations, and the right-hand side vector ${\bm b}$ of the AVEs
(\ref{Ku1}) is taken in such a way that the vector ${\bm x}^{*} = (x_1, x_2,\ldots,
x_n)^H$ with
\begin{equation}
x_k = (-1)^k\iota,\quad\ \ k = 1,2,\ldots,n
\label{kux1x}
\end{equation}
is the exact solution. The two sub-systems of linear equations involved are solved in the way if
$A{\bm x}^{*} = {\bm b}$, then ${\bm x}^{*} = A^{-1}{\bm b}$. Moreover, if the two sub-systems of
linear equations involved in the Picard-CSCS and the nonlinear CSCS-like iteration methods are solved
by exploiting the method introduced in \cite{MKCO} and using parallel computing, numerical performances
of the Picard-CSCS and the nonlinear CSCS-like iteration methods should become better.

On the other hand, Mangasarian modified the classical Newton iteration method for solving
AVEs by introducing the auxiliary diagonal matrix $\hat{D}({\bm x}) = \partial|{\bm x}| =
\mathrm{diag}(\mathrm{sign}({\bm x}))$, refer to \cite{OLMA} for details; then he established
the generalized Newton iterative scheme with the initial guess ${\bm x}^{(0)}$,
\begin{equation}
{\bm x}^{(k + 1)} = (A - \hat{D}({\bm x}^{(k)}))^{-1}{\bm b},
\label{eqx2s}
\end{equation}
so it notes that we need to solve a system of linear equations with the coefficient matrix
$J^{(k)} = A - \hat{D}({\bm x}^{(k)})$, i.e., Eq. (\ref{eqx2s}). If the matrix $J^{(k)} =
A - \hat{D}({\bm x}^{(k)})$ is very sparse, then Eq. (\ref{eqx2s}) can be solved by using
MATLAB's function ``$\setminus$". If the matrix $J^{(k)} = A - \hat{D}({\bm x}^{(k)})$ is
large-scale (even dense), the Eq. (\ref{eqx2s}) can be solved by using Krylov subspace methods,
such as GMRES \cite{YSMHS} and TFQMR \cite{RWFAT}. This consideration just follows the recent
method named the inexact semi-smooth Newton method, which has been introduced in \cite{JYBCO}.
In our numerical experiments, we also give the compared results between the proposed method
and the above generalized Newton iterative scheme.

In practical implementations, the optimal parameter $\sigma_{ \mathrm{HSS} } = \sqrt{\lambda_{
\mathrm{\mathrm{max}}}\lambda_{\mathrm{\mathrm{min}}}}$ recommended in \cite{ZZBGHG2} is employed
for the Picard-HSS and nonlinear HSS-like iteration methods, where $\lambda_{\mathrm{min}}$ and $\lambda_{
\mathrm{max}}$ are the minimum and the maximum eigenvalues of the Hermitian part $H$ of the matrix
$A$. Similarly, we adopt the optimal parameter $\sigma_{\mathrm{CSCS}}$ given in \cite{MKNC,WQSLL}
for the Picard-CSCS iteration method and the nonlinear CSCS-like iteration method.
More precisely, in our calculations, $\sigma_{\mathrm{CSCS}}$ is chosen according to the following formula
\begin{equation*}
\sigma_{\mathrm{CSCS}} =
\begin{cases}
\sqrt{\gamma_{min}\gamma_{max} - \zeta^{2}_{max}},& \mathrm{for}~~\zeta_{max} <
\sqrt{\gamma_{min}\gamma_{max}},\\
\sqrt{\gamma^{2}_{min} + \zeta^{2}_{max}},& \mathrm{for}~~\zeta_{max} \geq
\sqrt{\gamma_{min}\gamma_{max}},
\end{cases}
\end{equation*}
where $\gamma_{min}$ and $\gamma_{max}$ are the lower and the upper bounds of the real part of
the eigenvalues of the matrices $C$ and $S$, and $\zeta_{max}$ is the upper bound of the absolute
values of the imaginary part of the eigenvalues of the matrices $C$ and $S$. Meanwhile, it should mention that two optimal parameters
$\sigma_{\mathrm{HSS}}$ and $ \sigma_{\mathrm{CSCS}}$ only minimize the bounds of the convergence
factors (not the spectral radiuses selves) of the HSS and CSCS iteration matrices, respectively \cite{MZZGFZ1}. Admittedly,
the optimal parameters are crucial for guaranteeing fast convergence of these parameter-dependent
iteration methods, but they are generally difficult to be determined, see e.g. \cite{ZZBGHG2,ZZBXY,MZZGFZ1,DKST}
for a discussion of these issues.

To show that the proposed iteration methods can also be efficiently applied to deal with the complex system
of AVEs (\ref{Ku1}), we first construct and test the following example, which is a system of AVEs with complex
Toeplitz matrix.
\vspace{1mm}

\noindent\textbf{Example 1}. We consider that $A \in \mathbb{C}^{n\times n}$ is a complex
non-Hermitian, sparse and positive definite Toeplitz matrix with the following form
\begin{equation}
A = \begin{pmatrix}
\gamma & c\iota & d\iota    \\
-1 - c\iota & \gamma & c\iota & d\iota   \\
-1 - d\iota &-1 - c\iota & \gamma & c\iota & d\iota \\
&\ddots &\ddots &\ddots &\ddots &\ddots\\
&&-1 - d\iota &-1 - c\iota & \gamma & c\iota & d\iota \\
&&&-1 - d\iota &-1 - c\iota & \gamma & c\iota \\
&&&&-1 - d\iota &-1 - c\iota & \gamma \\
\end{pmatrix},
\label{matrix2}
\end{equation}
where $\iota = \sqrt{-1}$ and $c, d, \gamma\in\mathbb{R}$ are three given parameters. It means that the matrices $A$ in the targeted AVEs are defined
as Eq. (\ref{matrix2}). According to the performances of HSS-based methods, see
\cite{DKST,MZZGFZ,JJZT}, compared with other early established methods, we compare the
proposed CSCS-based methods with HSS-based methods in Example 1. Then we choose
different parameters $c$ and $d$ and present the corresponding numerical results in
Tables \ref{tab2}-\ref{tab3}.

\begin{table}[t]\small\tabcolsep=4.5pt
\begin{center}
\caption{The optimal parameters $\sigma^{*}_{opt}$ for Example 1.}
\begin{tabular}{ccccccccc}
\hline $\gamma$ &$(c,~d)$ &$\sigma^{*}_{opt}$       &$n = 128$ &$n = 256$ &$n = 512$ &$n = 1024$ &$n = 2048$ &$n = 4096$ \\
\hline
            10  &$(2,~3)$ &$\sigma_{\mathrm{HSS}}$  &2.9710    &2.9524    &2.9477    &2.9465     &2.9462     &2.9461     \\
                &         &$\sigma_{\mathrm{CSCS}}$ &1.1817    &1.1818    &1.1813    &1.1813     &1.1813     &1.1813     \\
            13.5 &$(3,~4)$ &$\sigma_{\mathrm{HSS}}$ &3.6871    &3.6595    &3.6525    &3.6507     &3.6503     &3.6502     \\
                &         &$\sigma_{\mathrm{CSCS}}$ &1.6008    &1.5997    &1.5989    &1.5989     &1.5988     &1.5989     \\
\hline
\end{tabular}
\label{tab1}
\end{center}
\end{table}

\begin{table}[!htbp]\small\tabcolsep=4.6pt
\begin{center}
\caption{Numerical results for Example 1 with order $n$, $\gamma = 10$, and $(c,d) = (2,3)$.}
\begin{tabular}{llllllll}
\hline Method&        &$n = 128$ &$n = 256$ &$n = 512$ &$n = 1024$ &$n = 2048$ &$n = 4096$  \\
\hline
Picard-HSS  &IT\_out &6            &6            &6            &6            &5            &5      \\
            &IT\_inn &9.8333       &9.8333       &9.5000       &9.5000       &9.2000       &9.2000  \\
            &IT      &59           &59           &57           &57           &46           &46      \\
            &CPU     &0.0158       &0.0215       &0.0243       &0.0338       &0.0511       &0.1044  \\
Picard-CSCS &IT\_out &6            &6            &6            &6            &6            &5       \\
            &IT\_inn &6.3333       &6.3333       &6.0000       &6.0000       &6.0000       &5.6000  \\
            &IT      &38           &38           &36           &36           &36           &28       \\
            &CPU     &0.0112       &0.0146       &0.0189       &0.0258       &0.0431       &0.0546   \\
HSS-like    &IT      &37           &36           &35           &34           &33           &31      \\
            &CPU     &0.0136       &0.0173       &0.0218       &0.0275       &0.0531       &0.1013   \\
CSCS-like   &IT      &24           &23           &22           &22           &21           &21       \\
            &CPU     &0.0072       &0.0097       &0.0131       &0.0215       &0.0276       &0.0514      \\
GN          &IT      &\texttt{Fail} &\texttt{Fail} &\texttt{Fail} &\texttt{Fail} &\texttt{Fail} &\texttt{Fail}      \\
            &CPU     &--           &--           &--           &--           &--           &--        \\
\hline
\end{tabular}
\label{tab2}
\end{center}
\end{table}

\begin{table}[!htbp]\small\tabcolsep=4.6pt
\begin{center}
\caption{Numerical results for Example 1 with order $n$, $\gamma = 13.5$, and $(c,d) = (3,4)$.}
\begin{tabular}{llllllll}
\hline Method &        &$n = 128$ &$n = 256$ &$n = 512$ &$n = 1024$ &$n = 2048$ &$n = 4096$ \\
\hline
Picard-HSS    &IT\_out &5         &5         &5         &5          &5          &5          \\
              &IT\_inn &10.0000   &10.0000   &10.0000   &10.0000    &10.0000    &10.0000    \\
              &IT      &50        &50        &50        &50         &50         &50         \\
              &CPU     &0.0129    &0.0161    &0.0206    &0.0298     &0.0544     &0.1135     \\
Picard-CSCS   &IT\_out &5         &5         &5         &5          &5          &5          \\
              &IT\_inn &6.2000    &6.0000    &6.0000    &6.0000     &5.8000     &5.8000     \\
              &IT      &31        &30        &30        &30         &29         &29         \\
              &CPU     &0.0098    &0.0123    &0.0166    &0.0248     &0.0332     &0.0581     \\
HSS-like      &IT      &41        &40        &39        &38         &37         &36         \\
              &CPU     &0.0121    &0.0174    &0.0192    &0.0306     &0.0587     &0.1131     \\
CSCS-like     &IT      &24        &23        &23        &22         &21         &21         \\
              &CPU     &0.0075       &0.0093       &0.0136    &0.0212     &0.0289     &0.0496     \\
GN            &IT      &\texttt{Fail} &\texttt{Fail} &\texttt{Fail} &\texttt{Fail} &\texttt{Fail} &\texttt{Fail}       \\
              &CPU     &--           &--           &--           &--           &--           &--   \\
\hline
\end{tabular}
\label{tab3}
\end{center}
\end{table}

Firstly, the optimal parameters $\sigma_{\mathrm{CSCS}}$ and $\sigma_{\mathrm{HSS}}$
for Example 1 are listed in Table \ref{tab1}. It is worth mentioning that with the
increase of the matrix dimension $n$, the optimal parameters $ \sigma_{ \mathrm{CSCS}} $ and
$\sigma_{\mathrm{HSS}}$ are almost fixed or decreasing slightly. Moreover, in Tables
\ref{tab2}-\ref{tab3}, we report numerical results with respect to the Picard-HSS,
the nonlinear HSS-like, the Picard-CSCS, the nonlinear CSCS-like iterations, and the generalized
Newton iterations using the MATLAB's function ``$\setminus$" (referred to as GN).
We also present the elapsed CPU time in seconds (denoted as CPU) and the number of
outer, inner and total iteration steps (outer and inner iterations only for both Picard-HSS and
Picard-CSCS) for the convergence performances (denoted as IT\_out, IT\_inn and IT,
respectively).

As seen from Tables \ref{tab2}-\ref{tab3}, it finds that except the GN method, the Picard-HSS, the
nonlinear HSS-like, the Picard-CSCS and the nonlinear CSCS-like iterative methods can successfully
achieve approximate solutions of the AVEs with all different matrix dimensions. When the dimension
$n$ is increasing, the number of outer and inner iteration steps are almost fixed for all iteration
methods, and the number of total iteration steps shows the similar phenomena. But the total CPU time
for all iteration methods are increasing quickly. Moreover, in terms of outer iteration steps, the
Picard-HSS iteration method and the Picard-CSCS iteration method have almost the same results, but
the Picard-CSCS iteration method is better than the Picard-HSS iteration method in terms of inner
iteration steps. Then as a result, the Picard-CSCS iteration method is also more competitive than
the Picard-HSS iteration method in aspects of the elapsed CPU time.

\begin{figure}[t]
\centering
\includegraphics[width=2.47in,height=2.45in]{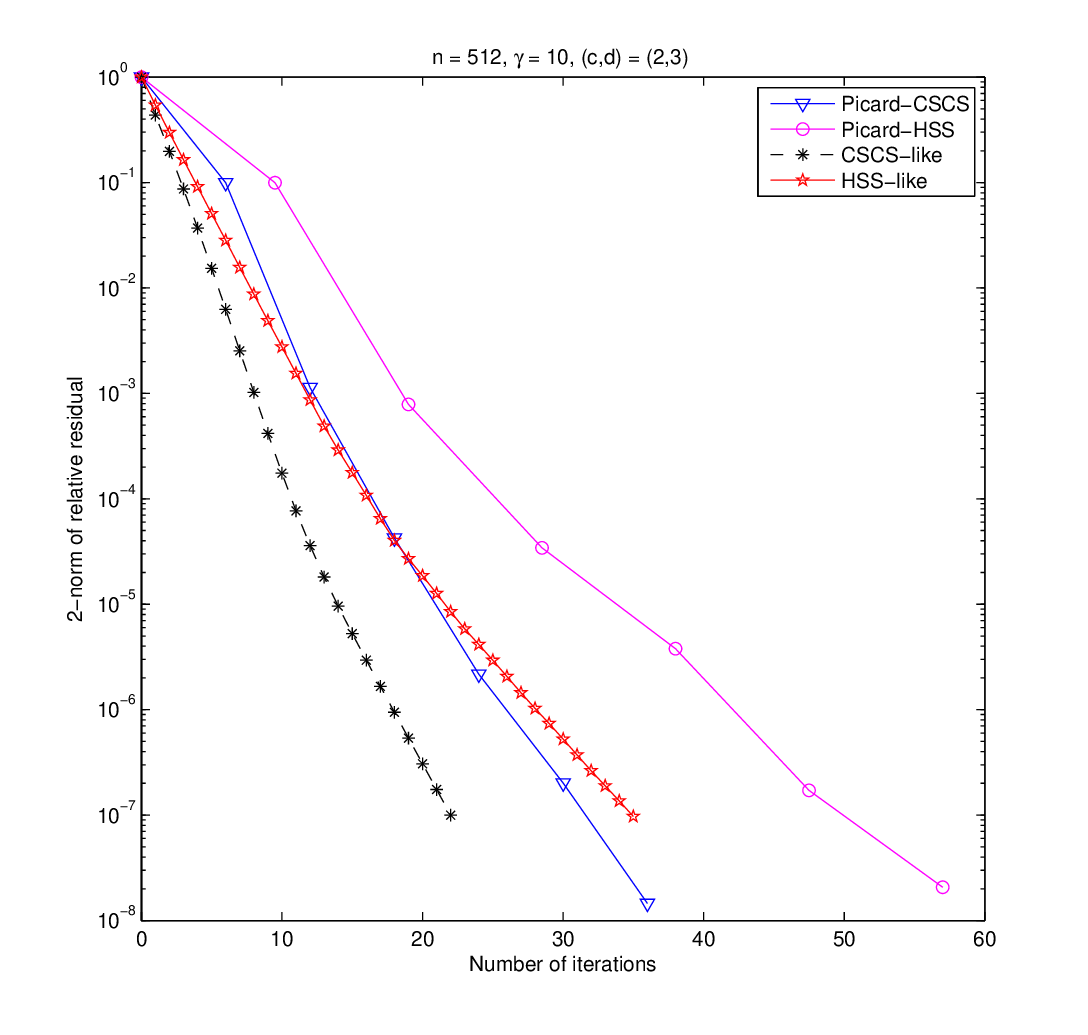}
\includegraphics[width=2.47in,height=2.45in]{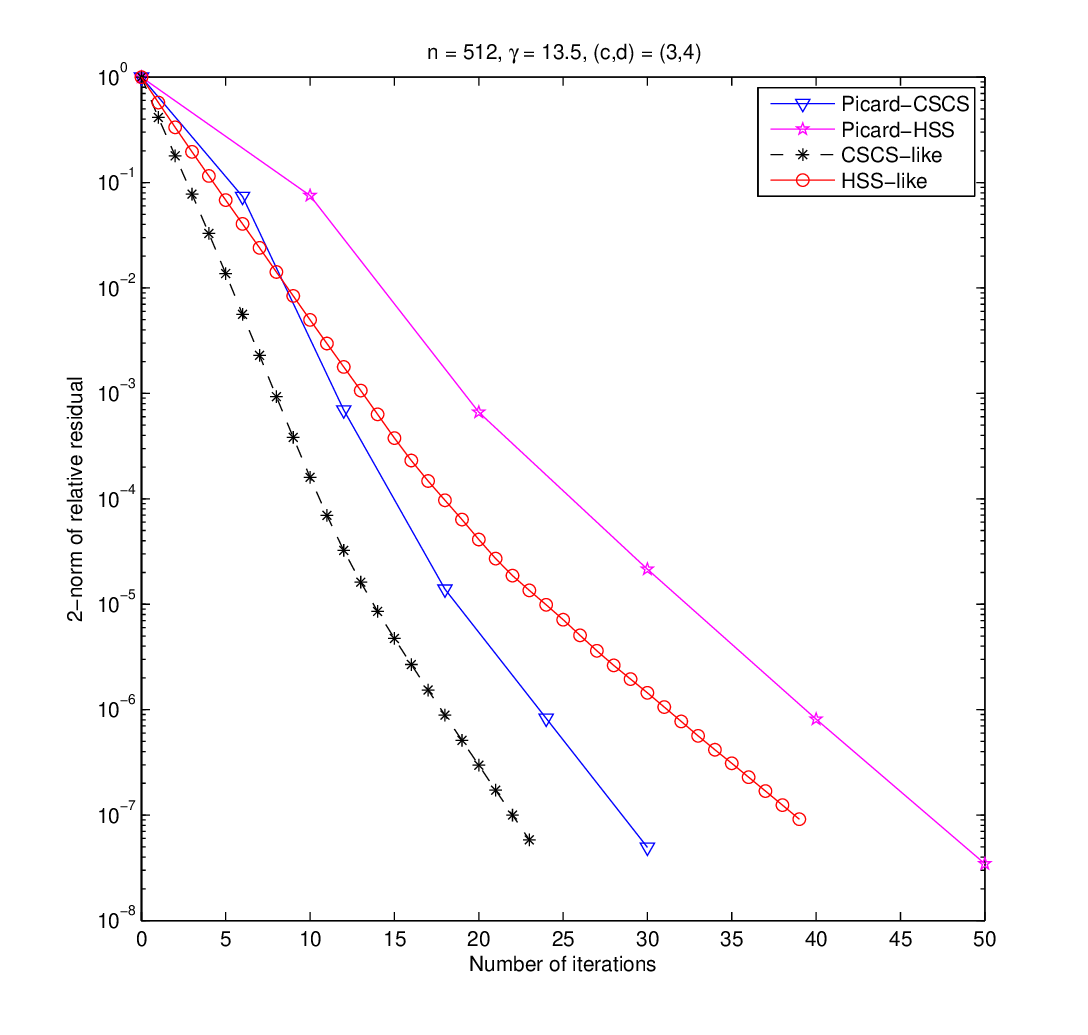}
\caption{{\small Convergence histories of the different iterative methods for two systems of
AVEs with the size $n = 512$ in Example 1.}}
\label{fig1x}
\end{figure}

On the other hand, from Tables \ref{tab2}-\ref{tab3}, we also observe that both the nonlinear
CSCS-like and the Picard-CSCS iteration methods are better than the nonlinear HSS-like and the
Picard-HSS iteration methods in terms of the number of iteration steps and the elapsed CPU time
for solving AVEs. In particular, the nonlinear CSCS-like method often enjoys the better performance
than the Picard-CSCS method in our implementations. Moreover, the convergence histories of residual
2-norms of these four different iterative algorithms are displayed in Fig. \ref{fig1x}, and the
performance profile based on CPU time for Example 1 with increasing the
matrix size $n$ is illustrated in Fig. \ref{fig2x}. In conclusion, the nonlinear CSCS-like
iteration method is the best choice for coping with AVEs concerning in Example 1. Besides,
the Picard-CSCS iteration method can be regarded as an acceptable alternative.

{\color{blue}\begin{figure}[!htbp]
\centering
\includegraphics[width=2.47in,height=2.4in]{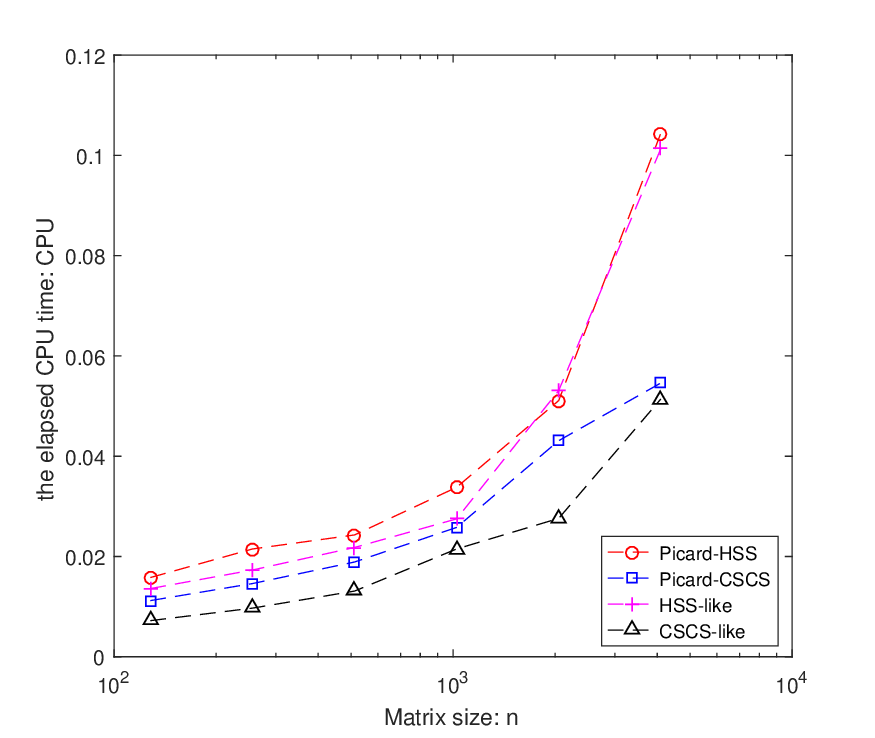}
\includegraphics[width=2.47in,height=2.4in]{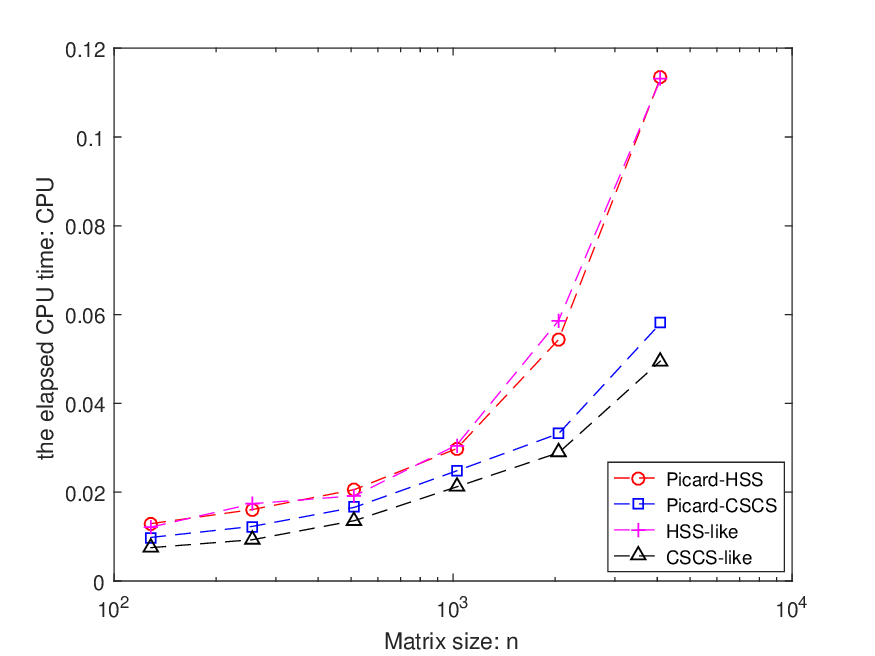}
\caption{{\small Performance profile based on CPU time under the matrix size $n$ in Example 1.}}
\label{fig2x}
\end{figure}}

\noindent\textbf{Example 2}. In order to evaluate the performances of the propose methods
comprehensively, we consider a family of the practical problems about the AVEs arising in
numerical solutions of the following one-dimensional nonlinear space fractional diffusion
equation, which is specially modified from Refs. \cite{HWCSKC,MMMCT},
\begin{equation}
\begin{cases}
\frac{\partial u(x,t)}{\partial t} = d_{+} \frac{\partial^{\alpha} u(x,t)}{\partial_{+} x^{\alpha}}
+ d_{-}\frac{\partial^{\alpha} u(x,t)}{\partial_{-} x^{\alpha}} + |u(x,t)|/\varsigma,\quad x\in (0,1),
\quad t\in [0,1],\\
u(0,t)= u(1,t)=0,\quad\ 0\leq t\leq 1,\\
u(x,0)= \phi(x),\quad\ \quad\ 0\leq x\leq 1,
\end{cases}
\label{eqr}
\end{equation}
where $\alpha\in (1,2)$ is the order of the fractional derivative, $\varsigma > 0$, and diffusion coefficients
$d_{\pm}$ are nonnegative; i.e., $d_{\pm} \geq 0$. Moreover, $\phi(x)$ is a known function. To solve Eq. (\ref{eqr})
numerically, let $N$ and $M$ be positive integers, and $h = 1/(N + 1)$ and $\tau = 1/M$ be the sizes of spatial
grid and time step, respectively. We define a spatial and temporal partition $x_j = jh$ for $j = 0,1,\ldots,N + 1$
and $t_m = m\tau$ for $m = 0,1,\ldots,M$. Let $u^{(m)}_j = u(x_j,t_m)$. In \cite{MMMCT}, Meerschaert and Tadjeran
proposed the shifted Gr\"{u}nwald approximation as follows,
\begin{subequations}
  \begin{align}
    \frac{\partial^{\alpha} u(x_j, t_m)}{\partial_{+} x^{\alpha}} & = \frac{1}{h^{\alpha}} \sum^{j +
    1}_{k = 0}g^{(\alpha)}_k u^{(m)}_{j - k + 1} + \mathcal{O}(h),\label{eqr2}\\
    \frac{\partial^{\alpha} u(x_j, t_m)}{\partial_{-} x^{\alpha}} & = \frac{1}{h^{\alpha}} \sum^{N
    - j + 2}_{k = 0}g^{(\alpha)}_k u^{(m)}_{j + k - 1} + \mathcal{O}(h),\label{eqr3}
  \end{align}
  \label{straincomponent}
\end{subequations}
where the coefficients $g^{(\alpha)}_k$ and corresponding properties are given in \cite[Proposition 1]{MMMCT,HKPHWS}.
Combining the implicit Euler scheme with Eqs. (\ref{straincomponent}) to discrete Eq. (\ref{eqr}),
then the final numerical scheme is
\begin{equation}
\frac{u^{(m)}_j - u^{(m - 1)}_j}{\tau} = \frac{d_{+}}{h^{\alpha}} \sum^{j + 1}_{k = 0}g^{(\alpha)}_k u^{(m)
}_{j - k + 1} + \frac{d_{-}}{h^{\alpha}} \sum^{N - j + 2}_{k = 0}g^{(\alpha)}_k u^{(m)}_{j + k - 1} +
|u^{(m)}_j|/\varsigma.
\label{eq15x}
\end{equation}
By using the similar ways given in \cite{HWCSKC}, it is not difficult to prove that the numerical scheme
(\ref{eq15x}) is unconditionally stable, which we will not pursue here. Let ${\bm u}^{(m)} = (u^{(m)}_1,
u^{(m)}_2,\ldots,u^{(m)}_N)^{T},~m = 0,1,\ldots,M$ and $I_N$ be the identity matrix of order $N$. Then
the numerical scheme (\ref{eq15x}) at the first temporal level $m = 1$ can be written in the following
matrix form
\begin{equation}
\Big[I_N - \frac{\tau}{h^{\alpha}}(d_{+}G_{\alpha} + d_{-}G^{T}_{\alpha})\Big]{\bm u}^{(1)} -
|{\bm u}^{(1)}| = {\bm u}^{(0)},
\end{equation}
where we take $\varsigma = \tau$ and $G_{\alpha}\in \mathbb{R}^{N\times N}$ is a nonsymmetric Toeplitz
matrix defined in \cite{HKPHWS}. According to Eq. (\ref{eq15x}), it implies that we need to handle a
system of nonlinear equations like the AVEs in (\ref{Ku1}) at each time step, i.e., there is a need for
solving the AVEs with the form $A{\bm u} - |{\bm u}| = {\bm u}^{(0)}$, where $A = I_N - \frac{\tau}{h^{
\alpha}}(d_{+}G_{\alpha} + d_{-}G^{T}_{\alpha})$ is also a nonsymmetric Toeplitz matrix. Meanwhile, for
simplicity, the vector ${\bm u}^{(0)}$ is still chosen as the same as that ${\bm x}^{*}$ in Eq. (\ref{kux1x})
is the solution of AVEs in (\ref{Ku1}).

Next, for the first temporal level $m = 1$, we employ these two CSCS-based iteration methods to solve
the above resultant AVEs, then the necessary condition for analyzing the convergence of the CSCS-based
iteration method is that both the circulant part $C$ and the skew-circulant part $S$ of $A$ are positive
definite. In fact, we have already mentioned that both the circulant part $C$ and the skew-circulant
part $S$ of the matrix $A = I_N - \frac{\tau}{h^{\alpha}}(d_{+}G_{\alpha} + d_{-}G^{T}_{\alpha})$ are
positive definite (see \cite{XMGTZH2} for details) via the similarly analyzed methods in \cite{WQSLL}.
It means that exploiting the CSCS-based iteration methods for solving the resulting AVEs is reasonable.
\textit{At the same time, it is worth mentioning that HSS-based iteration methods are not suitable for
Example 2 due to the Toeplitz coefficient matrix. Otherwise, it will lead to the complicated computations
for solving two sub-systems with the dense coefficient matrices $\sigma I + H$ and $\sigma I + S$}. In this example, since the matrix $J^{(k)}$ is a Toeplitz-plus-diagnoal
matrix, so there are no fast direct solvers for $J^{(k)}\tilde{{\bm x}} = {\bm u}^{(0)}$\footnote{It
is mainly because the displacement rank of the matrix $J^{(k)}$ can take any value between $0$ and $n$.
Hence, fast Toeplitz direct solvers that are based on small displacement rank of matrices cannot be
applied \cite[p.142]{MKNI}.}. Fortunately, it should note that the matrix-vector product
involving $J^{(k)}$ can be implemented via FFTs due to having the Toeplitz part $A$. It tells
us that the Krylov subspace methods can be compatibly exploited for solving $J^{(k)}\tilde{{\bm x}}
= {\bm u}^{(0)}$ at each iteration step, we denote them as the GN-TFQMR method and the
GN-GMRES method. In conclusion, we will compared the proposed CSCS-based iteration method with both
the GN-GMRES and GN-TFQMR methods for solving the resultant AVEs in Example 2. Numerical results are
reported in the following tables under different values of $\alpha,d_{\pm}$ and $h = \tau$. The total
number of (inner) iteration steps used for both GMRES and TFQMR methods, which are used to solve $J^{(k)}
\tilde{{\bm x}} = {\bm u}^{(0)}$, is no more than 15 in our practical implementations.

\begin{table}[!htbp]\small\tabcolsep=4.6pt
\begin{center}
\caption{The optimal parameters $\sigma^{*}_{opt}$ of the CSCS iteration method in Example 2.}
\begin{tabular}{cccccccc}
\hline $\alpha$ &$(d_{+},~d_{-})$  &\multicolumn{5}{c}{$\sigma^{*}_{opt}$}\\
[-2pt]\cmidrule(r{0.5em}){3-8} \\[-11pt]
    &          &$N = 128$ &$N = 256$ &$N = 512$ &$N = 1024$ &$N = 2048$ &$N = 4096$ \\
\hline
1.2 &$(0.5,~0.8)$  &1.4499    &1.5338    &1.6233    &1.7180     &1.8175 & 1.9216      \\
1.5 &$(0.6,~0.4)$ &2.7848    &3.2426    &3.7564    &4.3094     &4.8598  & 5.3131    \\
1.8 &$(0.7,~0.3)$ &5.8492    &6.9416    &7.0941    &15.3896    &26.7749 & 46.6033    \\
\hline
\end{tabular}
\label{tab4}
\end{center}
\end{table}

First of all, the optimal parameters $\sigma_{\mathrm{CSCS}}$ for Example 2 are
listed in Table \ref{tab4}. It is remarked that with the increase of the matrix
dimension $n$, the optimal parameters $\sigma_{\mathrm{CSCS}}$ are almost fixed
or increasing slightly for the cases of $\alpha = 1.2$ and $\alpha = 1.5$. Since
the case of $\alpha = 1.8$ corresponding to the coefficient matrix $A$ is very
ill-conditioned, so the optimal parameters $\sigma_{\mathrm{CSCS}}$ are varied
intensely. Moreover, in Tables \ref{tab5}-\ref{tab7}, we report the numerical
results with respect to the Picard-CSCS, nonlinear CSCS-like, GN-GMRES and GN-TFQMR
iterative methods. Similar to Example 1, we report the elapsed CPU time in seconds
and the number of outer, inner and total iteration steps (outer and inner iterations
only for Picard-CSCS, GN-GMRES and GN-TFQMR) for showing the convergence performances.

\begin{table}[!htbp]\small\tabcolsep=4.6pt
\begin{center}
\caption{Numerical results for Example 2 with order $N$, $\alpha = 1.2$, and $(c,d) = (0.5,0.8)$.}
\begin{tabular}{llllllll}
\hline Method &        &$N = 128$     &$N = 256$     &$N = 512$     &$N = 1024$    &$N = 2048$    &$N = 4096$    \\
\hline
Picard-CSCS   &IT\_out &6             &6             &6             &6             &6             &6             \\
              &IT\_inn &4.0000        &4.0000        &4.0000        &4.1667        &5.0000        &5.0000        \\
              &IT      &24            &24            &24            &25            &30            &30            \\
              &CPU     &0.0077        &0.0098        &0.0136        &0.01998       &0.0331        &0.0585        \\
CSCS-like     &IT      &12            &13            &14            &15            &16            &18            \\
              &CPU     &0.0032        &0.0051        &0.0078        &0.0099        &0.01776       &0.0334        \\
GN-GMRES      &IT\_out &\texttt{max}  &\texttt{max}  &\texttt{max}  &\texttt{max}  &\texttt{max}  &\texttt{max}  \\
              &IT\_inn &--            &--            &--            &--            &--            &--            \\
              &IT      &\texttt{Fail} &\texttt{Fail} &\texttt{Fail} &\texttt{Fail} &\texttt{Fail} &\texttt{Fail} \\
              &CPU     & --           &--            &--            &--            &--            &--            \\
GN-TFQMR      &IT\_out &\texttt{max}  &\texttt{max}  &\texttt{max}  &\texttt{max}  &\texttt{max}  &\texttt{max}  \\
              &IT\_inn &--            &--            &--            &--            &--            &--            \\
              &IT      &\texttt{Fail} &\texttt{Fail} &\texttt{Fail} &\texttt{Fail} &\texttt{Fail} &\texttt{Fail} \\
              &CPU     &--            &--            &--            &--            &--            &--            \\
\hline
\end{tabular}
\label{tab5}
\end{center}
\end{table}

\begin{table}[!htbp]\small\tabcolsep=4.6pt
\begin{center}
\caption{Numerical results for Example 2 with order $N$, $\alpha = 1.5$, and $(d_{+},d_{-}) = (0.6,0.4)$.}
\begin{tabular}{llllllll}
\hline Method &        &$N = 128$     &$N = 256$     &$N = 512$     &$N = 1024$    &$N = 2048$    &$N = 4096$    \\
\hline
Picard-CSCS   &IT\_out &6             &6             &6             &6             &6             &6             \\
              &IT\_inn &7.0000        &8.1667        &9.1667        &10.6667       &13.0000       &14.3333       \\
              &IT      &42            &49            &55            &64            &78            &86            \\
              &CPU     &0.0161        &0.0189        &0.0234        &0.0335        &0.0609        &0.1307        \\
CSCS-like     &IT      &24            &29            &35            &43            &54            &69            \\
              &CPU     &0.0063        &0.0084        &0.0119        &0.0218        &0.0486        &0.1133        \\
GN-GMRES   &IT\_out &\texttt{max}  &\texttt{max}  &\texttt{max}  &\texttt{max}  &\texttt{max}  &\texttt{max}  \\
              &IT\_inn &--            &--            &--            &--            &--            &--            \\
              &IT      &\texttt{Fail} &\texttt{Fail} &\texttt{Fail} &\texttt{Fail} &\texttt{Fail} &\texttt{Fail} \\
              &CPU     & --           &--            &--            &--            &--            &--            \\
GN-TFQMR      &IT\_out &\texttt{max}  &\texttt{max}  &\texttt{max}  &\texttt{max}  &\texttt{max}  &\texttt{max}  \\
              &IT\_inn &--            &--            &--            &--            &--            &--            \\
              &IT      &\texttt{Fail} &\texttt{Fail} &\texttt{Fail} &\texttt{Fail} &\texttt{Fail} &\texttt{Fail} \\
              &CPU     &--            &--            &--            &--            &--            &--            \\
\hline
\end{tabular}
\label{tab6}
\end{center}
\end{table}

\begin{table}[!htbp]\small\tabcolsep=4.6pt
\begin{center}
\caption{Numerical results for Example 2 with order $N$, $\alpha = 1.8$, and $(d_{+},d_{-}) = (0.7,0.3)$.}
\begin{tabular}{llllllll}
\hline Method &        &$N = 128$     &$N = 256$     &$N = 512$     &$N = 1024$    &$N = 2048$    &$N = 4096$    \\
\hline
Picard-CSCS   &IT\_out &5             &6             &10            &8             &9             &13            \\
              &IT\_inn &14.4000       &15.0000       &15.0000       &15.0000       &15.0000       &15.0000       \\
              &IT      &72            &90            &150           &120           &135           &195           \\
              &CPU     &0.0185        &0.0252        &0.0415        &0.0505        &0.0987        &0.2817        \\
CSCS-like     &IT      &59            &86            &146           &117           &117           &118           \\
              &CPU     &0.0101        &0.0179        &0.0386        &0.0517        &0.1007        &0.1911        \\
GN-GMRES      &IT\_out &\texttt{max}  &\texttt{max}  &\texttt{max}  &\texttt{max}  &\texttt{max}  &\texttt{max}  \\
              &IT\_inn &--            &--            &--            &--            &--            &--            \\
              &IT      &\texttt{Fail} &\texttt{Fail} &\texttt{Fail} &\texttt{Fail} &\texttt{Fail} &\texttt{Fail} \\
              &CPU     & --           &--            &--            &--            &--            &--            \\
GN-TFQMR      &IT\_out &\texttt{max}  &\texttt{max}  &\texttt{max}  &\texttt{max}  &\texttt{max}  &\texttt{max}  \\
              &IT\_inn &--            &--            &--            &--            &--            &--            \\
              &IT      &\texttt{Fail} &\texttt{Fail} &\texttt{Fail} &\texttt{Fail} &\texttt{Fail} &\texttt{Fail} \\
              &CPU     &--            &--            &--            &--            &--            &--            \\
\hline
\end{tabular}
\label{tab7}
\end{center}
\end{table}

Based on numerical results in Tables \ref{tab5}-\ref{tab7}, it finds that these two iterative solvers,
i.e., the Picard-CSCS and the nonlinear CSCS-like, can successfully obtain approximate solutions to the
AVEs for all different matrix dimensions; whereas both the GN-GMRES and GN-TFQMR iterative methods fully
fail to converge. It is mainly because the Newton-like iterative methods are usually sensitive to the
initial guess and the accuracy of solving the inner linear system corresponding to (\ref{eqx2s}) per
iterative step. When the matrix dimension $N$ is increasing, the number of outer iteration steps are almost
fixed or increasing slightly for all iteration methods, whereas the number of inner iteration steps show
the contrary phenomena for the cases with $\alpha = 1.2$ and $\alpha = 1.5$. Meanwhile, the total CPU time
and the total iteration steps for both the Picard-CSCS and the nonlinear CSCS-like iteration methods are
increasing quickly except the cases of $\alpha = 1.8$ with $N = 1024$ and $N = 2048$. On the other hand,
from Tables \ref{tab5}-\ref{tab7}, we also observe that the nonlinear CSCS-like method is almost more
competitive than the Picard-CSCS iteration method in terms of the number of iterations and the elapsed CPU
time for solving the AVEs. In particular, we can find that the nonlinear CSCS-like iteration method can require
slightly less number of iterations to converge than the Picard-CSCS iterative solver, but the Picard-CSCS
iterative solver can save a little elapsed CPU time with compared to the nonlinear CSCS-like iteration method
in our implementations. However, it still concludes that the nonlinear CSCS-like iterative method is the first
choice for solving the AVEs concerning in Example 2. At the same time, the Picard-CSCS iteration method can
be considered as a possible alternative.
\section{Conclusions}
\label{sec5}
In this paper, we have constructed two CSCS-based iteration methods for solving AVEs (\ref{Ku1})
with non-Hermitian Toeplitz matrix. Two CSCS-based iteration methods are based on separable
property of the linear term $A{\bm x}$ and the nonlinear term $|{\bm x}| + {\bm b}$ as well as on 
the CSCS of the involved non-Hermitian positive definite Toeplitz matrix $A$. By leveraging the
theory of nonsmooth analysis, the local convergence of nonlinear CSCS-like
iteration method has been investigated. Further numerical experiments have
shown that the Picard-CSCS and nonlinear CSCS-like iteration methods are feasible
and efficient nonlinear solvers for the AVEs. In particular, the nonlinear CSCS-like
iteration method often does better than the Picard-CSCS iteration method for solving
the AVEs. Finally, it is worth mentioning that how to employ suitable acceleration techniques
\cite{ZZBGHG1,JJZT,HFWPN} for enhancing the convergence of CSCS-based iteration methods, which
are affiliated with the fixed-point iteration, can remain an interesting topic of further research.
\section*{Acknowledgements}
We are grateful to the anonymous referees and editors for their insightful suggestions
and comments that improved the presentation of this paper.

\end{document}